\documentclass[12pt,reqno]{amsart}

\newcommand{\N}{{\mathbb N}}
\newcommand{\R}{{\mathbb R}}
\newcommand{\Z}{{\mathbb Z}}

\newcommand{\D}{\triangle}

\numberwithin{equation}{section}
\newtheorem{theorem}{Theorem}[section]
\newtheorem{defn}[theorem]{Definition}
\newtheorem{lemma}[theorem]{Lemma}
\newtheorem{remark}[theorem]{Remark}

\begin{document}

\title[Invariant manifold reduction and dynamical approximation ]
{A dynamical approximation for stochastic partial differential
equations}

\author[W.  Wang \& J. Duan ]
{Wei  Wang \& Jinqiao Duan    }

\address[W.~Wang ]
{Institute of Applied Mathematics\\Chinese Academy of Sciences\\
Beijing, 100080, China; Department of Mathematics \\ Nanjing
University\\ Jiangsu, 210093, China}
\email[W.~Wang]{wangweinju@yahoo.com.cn}

\address[J.~Duan]
{Department of Applied Mathematics\\
Illinois Institute of Technology\\
Chicago, IL 60616, USA} \email[J.~Duan]{duan@iit.edu}

 \date{July 16, 2007(revised version); May 4, 2007(original version)}

\thanks{The authors would like to thank Dirk Bl\"omker, Tomas Caraballo
and Peter E. Kloeden for helpful comments. This work was partly
supported by the NSF Grants DMS-0209326, DMS-0542450, NSFC Grant
No-10626052 and the Outstanding Overseas Chinese Scholars Fund of
the Chinese Academy of Sciences.}



\begin{abstract}
Random invariant manifolds often provide geometric structures for
understanding stochastic dynamics. In this paper, a dynamical
approximation estimate is derived for a class of stochastic
 partial differential equations, by showing that the random invariant manifold
is   almost surely asymptotically complete.  The asymptotic
dynamical behavior is thus described by a stochastic
 ordinary differential system on the random invariant manifold,
 under suitable conditions. As an application, stationary
 states (invariant measures) is considered for one example of
    stochastic   partial differential equations.
\end{abstract}

\maketitle

\section{Introduction}\label{s1}

Stochastic partial differential equations (SPDEs or stochastic
PDEs) arise as macroscopic mathematical models of complex systems
under random influences. There have been rapid progresses in this
area \cite{Gar, Roz, PZ92, WaymireDuan, B05, HuangYan}. More
recently, SPDEs have been investigated in the context of random
dynamical systems (RDS) \cite{Arn98}; see for example,
\cite{CLR00,CDSch03,Cheu95,CF94, CFD97, Sch92, DLSch03, DLSch04},
among others.


Invariant manifolds are special invariant sets represented by graphs
in  state spaces (function spaces) where solution processes of SPDEs
live. A random invariant manifold  provides a geometric structure to
reduce stochastic dynamics. Stochastic bifurcation, in a
  sense, is about the   changes in  invariant structures
for     random dynamical systems. This  includes   qualitative
 changes of   invariant manifolds, random attractors, and
    invariant  measures or stationary states.

Duan et  al. \cite{DLSch03, DLSch04} have recently proved results on
existence and smoothness of random invariant manifolds for a class
of stochastic partial differential equations. In this paper, we
further derive a dynamical approximation estimate between the
solutions of stochastic partial differential equations \emph{and}
the orbits on the random invariant manifolds. This is achieved by
showing that the random invariant manifold is   almost surely
asymptotically complete (see Definition \ref{completeness}).  The
asymptotic dynamical behavior thus can be described by a stochastic
ordinary differential system on the random invariant manifold, under
suitable conditions. In this approach one key assumption is  that
the global Lipschitz constant of nonlinear term is small enough.

If the invariant manifold is almost surely asymptotically complete
 we can approximate the infinite dimensional system by a
system restricted on the random invariant manifold which is in fact
   finite dimensional. That is the infinite dimensional
system is reduced to a finite dimensional system, which is useful
for understanding asymptotic behavior of the original stochastic
system \cite{roberts}.





\bigskip

As a application,  in \S 5, we study the existence of stationary
solutions of  a hyperbolic SPDE. Specifically, we will consider the
following stochastic hyperbolic
 partial differential  equation, in \S \ref{s5}, with large diffusivity and highly
damped term, on the space-time domain $[0, 2\pi]\times (0,
+\infty)$
\begin{equation}\label{Dir}
u_{tt}(t,x)+\alpha u_t(t,x)=\nu\Delta u(t, x)+f(u(t,x), x)+u(t,x)
\circ  \dot{W}(t),
\end{equation}
with
$$
u(0, x)=u_0,\;\;u_t(0, x)=u_1,\;\;u(t, 0)=u(t, 2\pi)=0, \;\; t>0
$$
where $\nu$ and $\alpha$ are both positive. And $f\in C^2(\R,\R)$ is
a bounded globally Lipschitz nonlinearity. Note that the stochastic
Sine-Gordon equation ($f=\sin u$) is an example. When the damping is
large enough the existence of the stationary solutions for
(\ref{Dir}) is obtained by considering the stochastic system on the
random invariant manifold.


\bigskip

   Invariant manifolds are often used as a
  tool to study the structure of attractors \cite{Constantin,Tem97}.
  In this paper we first consider invariant
manifolds   for a class infinite dimensional RDS defined by SPDEs,
   then reduce the random dynamics to the invariant
manifolds.  When the invariant manifolds are shown to be almost
surely asymptotically complete (see Definition
\ref{completeness}), we obtain dynamical approximations of the
solutions of stochastic PDEs by orbits on the invariant manifolds.
  Almost sure cone invariance concept (see Definition
\ref{cone}) is used to prove almost sure asymptotic completeness
property of the random invariant manifolds.


\bigskip

 This paper is organized as follows. We state the main result on
    dynamical approximation  for a class of
 SPDEs in section \ref{s2}.    Then we  recall   background materials
 in random dynamical systems and the existence result  of invariant manifolds \cite{DLSch03}
 in section \ref{s3}.   The main result is proved in section \ref{s4}, and
applications in    detecting stationary states are
discussed in the final section \ref{s5}.\\


\section{Main result}\label{s2}

We consider the stochastic evolutionary  system
\begin{equation}\label{e1}
du(t)=(Au(t)+F(u(t)))dt+  u(t)\circ dW(t)
\end{equation}
where $A$ is the generator of a $C_0$-semigroup $\{e^{tA}\}_{t\geq
0}$ on real valued separable Hilbert space $(H, |\cdot|)$ with inner
product $\langle,\cdot, \rangle$;   $F:\R\rightarrow \R$ is a
continuous nonlinear function with $F(0)=0$ and Lipschitz constant
$L_F$ is assumed to be small; and $W(t)$ is a standard real valued
Wiener process. Moreover, $\circ$ denotes the stochastic
differential in the sense of Stratonovich. Suppose that $\sigma(A)$,
the spectrum of operator $A$, splits as
\begin{equation}\label{spetrum}
 \sigma(A)=\{\lambda_k,\;k\in\N \}=\sigma_c \cup\sigma_s,\;\;\sigma_c, \sigma_s\neq  \emptyset
\end{equation} with
$$
\sigma_c\subset \{z\in \mathbb{C}:Re z\geq 0\},\;\; \sigma_s\subset
\{z\in\mathbb{C}:Re z<0\}
$$
where $\mathbb{C}$ denotes the complex numbers set. $\sigma_c$ is
assumed to be a finite set. Denote the corresponding eigenvectors to
$\{\lambda_k,\;k\in\N\}$ by $\{e_1, \cdots, e_n, e_{n+1}, \cdots\}$.
By the above assumptions, there is an $A$-invariant decomposition
$H=H_c\oplus H_s$ such that for the restrictions $A_c=A|_{H_c}$,
$A_s=A|_{H_s}$ one has $\sigma_c=\{z: z\in\sigma(A_c)\}$ and
$\sigma_s=\{z: z\in\sigma(A_s)\}$. Moreover $\{e^{tA_c}\}$ is a
group on $H_c$ and there exist  projections $\Pi_c$ and $\Pi_s$ such
that $\Pi_c+\Pi_s=Id_H$, $A_c=A|_{Im \Pi_c}$ and $A_s=A|_{Im
\Pi_s}$. We also suppose that there are
postive constants $\alpha$ and
$\beta$ with property $0\leq \alpha<\beta$, such that
\begin{equation}\label{expdic1}
|e^{tA_c}x|\leq  e^{\alpha t}|x|, \;\;t\leq 0,
\end{equation}
\begin{equation}\label{expdic2}
|e^{tA_s}x|\leq  e^{-\beta t}|x|,\;\; t\geq 0.
\end{equation}
For instance,   $-A$ may be a
strongly elliptic and symmetric second order differential operator
on a smooth domain with zero Dirichlet boundary condition.

Since  $\sigma_c$ is a finite set, $H_c$ is a finite dimensional
space of dimension, say, $\dim H_c=n$. For $u\in H$ we have
$u=u_c+u_s$ with $u_c=\Pi_c u\in H_c$ and $u_s=\Pi_s u\in H_s$.
Furthermore, we assume that the projections $\Pi_c$ and $\Pi_s$
commute with $A$. Define the nonlinear map on $H_c$ as
\begin{eqnarray*}
F_c&:& H_c\rightarrow H_c\\
&& u_c\mapsto F_c (u_c)=\sum_{i=1}^n\langle F(u_c+0), e_i \rangle
e_i
\end{eqnarray*}
where $0$ is the zero element in the vector space $H_s$. The concept
of random invariant manifolds will be introduced in the next
section. Let $\{\theta_t\}_{t\in\R}$ be the metric dynamical system
generated by Wiener process $W(t)$, see (\ref{MD}). We will obtain
the following main result.

\medskip

\begin{theorem} \label{theorem}
(\textbf{Dynamical approximations}) \\
 Consider the following stochastic evolutionary system
\begin{equation} \label{spde2}
du(t)= [Au(t)+F(u(t))]dt  +    u(t)\circ dW(t),
\end{equation}
where the linear operator $A$ and the nonlinearity $F$ satisfy the
conditions listed above. If  the Lipschitz constant of the
nonlinearity $F$  is small enough, then this stochastic system has
an n dimensional invariant manifold $\mathcal{M}(\omega)$ and there
exists a positive random variable $D(\omega)$ and a positive
constant $k$ such that: For any solution $u(t,\theta_{-t}\omega)$ of
(\ref{spde2}), there is an orbit   $U(t,\theta_{-t}\omega)$ on the
invariant manifold $\mathcal{M}( \omega)$, with the following
approximation property:
\begin{equation}\label{appro}
|u(t, \theta_{-t}\omega)-U(t,\theta_{-t}\omega)|\leq
D(\omega)|u(0)-U(0)|e^{-kt}, \;\; t>0, \;\;almost \;\;surely.
\end{equation}
\end{theorem}


\begin{remark}\label{remark}
The above result can be seen as a dynamical approximation for the
system (\ref{spde2}). Any solution $u$ of (\ref{spde2}) can be
approximated by an orbit $U$ on the manifold
$\mathcal{M}(\omega)$. In fact the function $U$  can be
represented as $u_c+\bar{h}^s(u_c)$, where $u_c$  satisfies the
following stochastic equation
\begin{equation}\label{e3}
du_c(t)=\left[A_c u_c(t)+F_c(u_c(t)+\bar{h}^s(u_c(t),
\theta_t\omega) ) \right]dt+ u_c(t) \circ dW(t),
\end{equation}
and, moreover,
\begin{equation*}
\bar{h}^s: H_c\rightarrow H_s
\end{equation*}
is a random Lipschitz map; see \S \ref{s4}. Here we also remark that
$\bar{h}^s$ depends on $\omega$, so (\ref{e3}) in fact is a
non-autonomous stochastic differential equation on $H_c$.

\end{remark}




\section{Random invariant manifolds}\label{s3}

Now we recall the basic concepts in random dynamical systems, and
the basic result on the existence of random invariant manifolds for
stochastic PDEs, from Duan et al. \cite{DLSch03, DLSch04}.


For our purpose we work on the canonical probability space
$(\Omega_0, \mathcal{F}_0,\mathbb{P})$, where the sample space
$\Omega_0$ consists of the sample paths of $W(t)$, that is
$$
\Omega_0=\{w\in C([0, \infty), \R): w(0)=0 \},
$$
for more see \cite{Arn98}. Consider the following stochastic
evolutionary equations
\begin{eqnarray} \label{SEE}
u_t&=&  Au+F(u)+  u \circ \dot{W}(t),\\
u(0)&=&u_0\in H.\nonumber
\end{eqnarray}
This equation can be written in the following mild integral form
\begin{equation}\label{mild}
u(t)=e^{ At}u_0+\int_0^t e^{A(t-s)}F(u(s))ds+\int_0^t
e^{A(t-s)}u(s)\circ dW(s)
\end{equation}
By the assumption of $A$ and $F$ we know that equation (\ref{mild})
has a unique solution $u(t,\omega;u_0)\in L^2(\Omega_0, C(0, T; H))$
for any $T>0$ in the sense of probability. For more about the
solution of SPDEs we refer to \cite{PZ92}.

\bigskip

    We now present some basics of random dynamical
systems. First we start with a driven dynamical system   which
models white noise:
$\theta_t:(\Omega_0,\mathcal{F}_0,\mathbb{P})\rightarrow
(\Omega_0, \mathcal{F}_0,\mathbb{P})$, $t\in\R$ that satisfies the
usual definition for a (deterministic) dynamical system
\begin{itemize}
    \item $\theta_0=id$,
    \item $\theta_t\theta_s=\theta_{t+s}$,\;\; for all $s$,
    $t\in\R$,
    \item the map $(t,\omega)\mapsto \theta_t\omega$ is
    measurable and $\theta_t\mathbb{P}=\mathbb{P}$ for all
    $t\in\R$.
\end{itemize}
 A random dynamical system (RDS) on a metric space $(\mathcal{X}; d)$ with
Borel $ $- algebra $\mathcal{B}$ over $\theta_t$ on $(\Omega_0,
\mathcal{F}_0, \mathbb{P})$ is a measurable map
\begin{eqnarray*}
   \varphi: \R^+\times\Omega_0\times \mathcal{X} &\rightarrow& \mathcal{X}\\
     (t,\omega,x)&\mapsto& \varphi(t,\omega)x
\end{eqnarray*}
 such that
\begin{itemize}
    \item[(i)] $\varphi(0,\omega)=id$ (on $\mathcal{X}$);
    \item[(ii)] $\varphi(t+s,\omega)=\varphi(t,\theta_s\omega)\varphi(s, \omega)$
    $\forall t,s \in\R^+$
    for almost all $\omega\in\Omega_0$ (cocycle property).
\end{itemize}
A RDS $\varphi$ is continuous or differentiable if $\varphi(t,
\omega) : \mathcal{X}\rightarrow \mathcal{X}$ is continuous or
differentiable (see \cite{Arn98} for more details on RDS). As we
consider the canonical probability space the driven dynamical
system $\theta_t$ can be defined as
\begin{equation}\label{MD}
\theta_t\omega(\cdot)=\omega(\cdot+t)-\omega(t),\;\; t\in \R,
\end{equation}
where $\omega(\cdot)\in \Omega_0$ is a sample path of the Wiener
process or Brownian motion $W(t)$.

For a continuous random dynamical system $\varphi(t, \omega) :
\mathcal{X}\rightarrow \mathcal{X}$ over
$\big(\Omega_0,\mathcal{F}_0,\mathbb{P},$ \\ $ (\theta_t)_{t\in
\R}\big)$, we  need the following  notions to describe its dynamical
behavior.

\begin{defn}
A collection $M=M(\omega)_{\omega\in\Omega}$, of nonempty closed
sets $M(\omega)$, $\omega\in\Omega$,  contained in $\mathcal{X}$,
is called a random  set if
\begin{equation*}
\omega\mapsto\inf_{y\in M(\omega)}d(x, y)
\end{equation*}
is a real valued random variable for any $x\in \mathcal{X}$.
\end{defn}

\begin{defn}
A random set $B(\omega)$ is called a tempered absorbing set for a
random dynamical system $\varphi$ if for any bounded set
$K\subset\mathcal{X}$ there exists $t_K(\omega)$ such that
$\forall t\geq t_K(\omega)$
\begin{equation*}
\varphi\big(t,\theta_{-t}\omega,K\big)\subset B(\omega).
\end{equation*}
and for all $\varepsilon>0$
\begin{eqnarray*}
 \lim_{t\rightarrow\infty}e^{-\varepsilon
t}d\big(B(\theta_{-t}\omega)\big)=0, \;\;a.e. \;\omega\in \Omega,
\end{eqnarray*}
where $d(B)=\sup_{x\in B}d(x, 0)$, with $0\in\mathcal{X}$, is the
diameter of $B$.
\end{defn}

For more about random set we refer to \cite{CF94}.

\begin{defn}
A random set $M(\omega)$ is called a positive invariant set for a
random dynamical system $\varphi(t,\omega,x)$ if
\begin{equation*}
\varphi(t,\omega, M(\omega))\subset
M(\theta_t\omega),\;\;for\;t\geq 0.
\end{equation*}
If $M(\omega)$ can be written as  a graph of a Lipschitz mapping
\begin{equation*}
\psi(\cdot,\omega): \mathcal{X}_1\rightarrow \mathcal{X}_2
\end{equation*}
where $\mathcal{X}=\mathcal{X}_1\oplus \mathcal{X}_2$, such that
\begin{equation*}
M(\omega)=\{x_1+\psi(x_1,\omega)|x_1\in \mathcal{X}_1 \}.
\end{equation*}
Then $M(\omega)$ is called a Lipschitz invariant manifold of
$\varphi$.
\end{defn}

For more about the random invariant manifold theory, see
\cite{DLSch03}.


In order to apply the random dynamical systems framework, we
transform the stochastic PDE (\ref{SEE}) into a random partial
differential equation (random PDE). To this end, we introduce the
following stationary process $z(\omega)$ which solves
\begin{equation}\label{OU}
dz+zdt=dW.
\end{equation}
Solving equation (\ref{OU}) with initial value
\begin{equation*}
z(\omega)=-\int_{-\infty}^0e^{\tau}w(\tau)d\tau,
\end{equation*}
 we have a unique stationary process for (\ref{OU})
\begin{equation*}
z(\theta_t\omega)=-\int_{-\infty}^0e^{\tau}\theta_{\tau}\omega(\tau)d\tau=
-\int_{-\infty}^0 e^{\tau}\omega(t+\tau)d\tau+\omega(t).
\end{equation*}
The mapping $t\mapsto z(\theta_t\omega)$ is continuous. Moreover
\begin{equation*}
\lim_{t\rightarrow\pm\infty}\frac{|z(\theta_t\omega)|}{|t|}=0\;\;{\rm
and }\;\; \lim_{t\rightarrow\pm\infty}\frac{1}{t}\int_0^t
z(\theta_\tau\omega)d\tau=0\;\; for\;a.e.\;\omega\in\Omega_0.
\end{equation*}
We should point out that above properties hold in a $\theta_t$
invariant set $\Omega\subset\Omega_0$ of full probability. For the
proof see \cite{DLSch03}. In the following part to the end we
consider the probability space $(\Omega, \mathcal{F}, \mathbb{P})$
where
$$
\mathcal{F}=\{\Omega\cap U: U\in \mathcal{F}_0\}.
$$

\bigskip

The following random PDE on the probability space $(\Omega,
\mathcal{F}, \mathbb{P})$ is a transformed version of the original
stochastic PDE  (\ref{SEE}).
\begin{equation}\label{rpde1}
v_t= Av+G(\theta_t\omega, v)+z(\theta_t\omega)v,\;\;v(0)=x\in H
\end{equation}
where $G(\omega, v):=e^{-z(\omega)}F(e^{z(\omega)}v)$. It is easy to
verify that $G$ has the same Lipschitz constant as $F$  and
$G(0)=0$. By the classical evolutionary equation theory equation
(\ref{rpde1}) has a unique solution $v(t,\omega; x)$ which is
continuous in $x$ for every $\omega\in\Omega$. Then
$$
(t, \omega, x)\mapsto \varphi(t,\omega)x:=v(t,\omega;x)
$$
defines a continuous random dynamical system.  We now introduce
the transform
\begin{equation} \label{transform}
T(\omega, x)=xe^{-z(\omega)}
\end{equation}
and its inverse transform
\begin{equation}
T^{-1}(\omega, x)=xe^{z(\omega)}
\end{equation}
for $x\in H$. Then for the random dynamical system $v(t,\omega;
x)$ generated by (\ref{rpde1}),
$$
(t,\omega, x)\rightarrow T^{-1}(\theta_t\omega, v(t,\omega;
T(\omega, x)):=u(t,\omega; x)
$$
is the random dynamical system generated by (\ref{SEE}). For more
about the relation between
(\ref{SEE}) and (\ref{rpde1}) we refer to \cite{DLSch03}.\\


\bigskip

We now prove the existence of a random invariant manifold for
(\ref{rpde1}) as in \cite{DLSch03, DLSch04}.

Let $d_H$ be the metric induced by the norm $|\cdot|$. Then $(H,
d_H)$ is a complete separable metric space.

Now we briefly give an approach to obtain the random invariant
manifold for system ({\ref{rpde1}}) by the Lyapunov-Perron method;
for detail see \cite{DLSch04}.

Projecting the system (\ref{rpde1}) onto $H_c$ and $H_s$
respectively we have
\begin{eqnarray}
\dot{v}_c&=&  A_c v_c+z(\theta_t\omega)v_c+ g_c(\theta_t\omega, v_c+v_s) \label{rpdec}\\
\dot{v}_s&=&  A_s v_s+z(\theta_t\omega)v_s+g_s(\theta_t\omega,
v_c+v_s) \label{rpdes}
\end{eqnarray}
where
$$
g_c(\theta_t\omega, v_c+v_s )=\Pi_c G(\theta_t\omega, v_c+v_s)
$$
and
$$
g_s(\theta_t\omega, v_c+v_s )=\Pi_sG(\theta_t\omega, v_c+v_s).
$$

 Define the following Banach space for $\eta$,
 $-\beta<\eta<-\alpha$,
\begin{eqnarray*}
C_\eta^-=\{v: (-\infty, 0]\rightarrow H:&v& \textrm{is continuous
and}\\ && \sup_{t\in(-\infty, 0]}e^{-\eta
t-\int_0^tz(\theta_s\omega)ds}|v(t)|< \infty \}
\end{eqnarray*}
with norm
$$
|v|_{C_{\eta}^-}=\sup_{t\in(-\infty,0]}e^{-\eta
t-\int_0^tz(\theta_s\omega)ds}|v(t)|.
$$

Define the nonlinear operator $\mathcal{N}$ on $C_\eta^-$ as
\begin{eqnarray}\label{T}
\mathcal{N}(v,\xi)(t, \omega)&=&e^{ A_ct+\int_0^tz(\theta_\tau
\omega)d\tau}\xi+\int_0^te^{
A_c(t-\tau)+\int_\tau^tz(\theta_\varsigma\omega)d\varsigma}g_c(\theta_\tau\omega,
v_c+v_s)d\tau \nonumber\\
&&+\int^t_{-\infty} e^{
A_s(t-\tau)+\int_\tau^tz(\theta_\varsigma\omega)d\varsigma}g_s(\theta_\tau\omega,
v_c+v_s)d\tau
\end{eqnarray}
where $\xi\in H_c$. Then for any given $\xi\in H_c$ and each $v$,
$\bar{v}\in C_\eta^-$, we have
\begin{eqnarray}
&&|\mathcal{N}(v,\xi)-\mathcal{N}(\bar{v},\xi)|_{C_\eta^-}\nonumber\\
&\leq&\sup_{t\leq 0}\Big\{ L_F\Big (\int_0^te^{(-\alpha-\eta)(t-s)}
ds+
 \int_\infty^te^{(-\beta-\eta)(t-s)}ds
\Big )\Big\} |v-\bar{v}|_{C_\eta^-}\nonumber\\
&\leq&
 L_F\big(\frac{1}{\eta+\beta}-\frac{1}{\alpha+\eta}\big)
|v-\bar{v}|_{C_\eta^-}.
\end{eqnarray}
If
\begin{equation}\label{c1}
 L_F\big(\frac{1}{\eta+\beta}-\frac{1}{\alpha+\eta}\big)<1,
\end{equation}
then by the fixed point argument
\begin{equation}\label{fix}
v=\mathcal{N}(v)
\end{equation}
has a unique solution $v^*(t,\omega;\xi)\in C_\eta^-$. Let
$h^s(\xi,\omega)=\Pi_sv^*(0, \omega; \xi)$. Then
\begin{equation}\label{h}
h^s(\xi, \omega)=\int^0_{-\infty} e^{
-A_s\tau+\int_\tau^0z(\theta_\varsigma\omega)d\varsigma}g_s(\theta_\tau\omega,
v^*(\tau, \omega;\xi))d\tau,
\end{equation}
 $h^s(0,\omega)=0$ and $h^s$ is Lipschitz continuous with Lipschitz constant
 $CL_F$, $C>0$ is a positive constant. Then we have the following result about the
existence of random invariant manifold for the random dynamical
system $\varphi(t,\omega)$ generated by (\ref{rpde1}). For the
detailed proof see \cite{DLSch04}.

\begin{lemma}\label{manifold}
Suppose the assumptions on $A$ and $F$ in section \ref{s2} and
condition (\ref{c1}) hold. Then there exists a Lipschitz continuous
random  invariant manifold $\mathcal{\widetilde{M}}(\omega)$ for
$\varphi(t,\omega)$ which is given by
$\mathcal{\widetilde{M}}(\omega)=\{\xi+h^s(\xi, \omega):\xi\in H_c
\}$.
\end{lemma}
Then by the transform $T$, as defined in (\ref{transform}), we
have the following conclusion.
 \begin{lemma} (\textbf{Invariant manifold for stochastic PDE})\\
$\mathcal{M}(\omega)=T^{-1}(\omega, \widetilde{M}(\omega))$ is a
Lipschitz continuous random invariant manifold for the stochastic
PDE (\ref{SEE}).
 \end{lemma}

\begin{remark}
Note that $\mathcal{M}(\omega)$ is independent of the
choice of $\eta$.\\
\end{remark}

\section{Dynamical approximations} \label{s4}

In this section we prove Theorem \ref{theorem}, by showing that
the invariant manifold obtained in the last section has the
\emph{almost sure  asymptotic completeness} property. Then the
dynamical behavior of (\ref{rpde1}) is determined by the system
restricted on the invariant manifold.


The following concept is important in the study of the dynamical
approximations of stochastic PDEs.

\begin{defn}\label{completeness}
(\textbf{Almost sure asymptotic completeness})\\
 Let $\mathcal{M}(\omega)$
be an invariant manifold for a random dynamical system
$\varphi(t,\omega)$. The invariant manifold $\mathcal{M}$ is
called almost surely asymptotically complete  if for every $ x\in
H$, there exists $y\in\mathcal{M}(\omega)$ such that
$$
|\varphi(t,\omega)x-\varphi(t,\omega)y|\leq D(\omega)|x-y|e^{-kt},
\;\;t\geq 0
$$
for almost all $\omega\in\Omega$, where $k$ is some positive
constant and $D$ is a positive   random variable.
\end{defn}

Now we introduce the \emph{almost sure cone invariance} concept.
For a positive random variable $\delta$, define the following
random set
\begin{equation*}
\mathcal{C}_\delta:=\big\{(v,\omega)\in H\times\Omega: |\Pi_s
v|\leq \delta(\omega) |\Pi_c v | \big \}.
\end{equation*}
And the fiber $\mathcal{C}_{\delta(\omega)}(\omega)=\{v:
(v,\omega)\in \mathcal{C}_\delta\}$ is called \textit{random
cone}. For a given random dynamical system $ \varphi(t,\omega)$ we
give the following definition.

\begin{defn}\label{cone}
(\textbf{Almost sure cone invariance})\\
 For a random cone
$\mathcal{C}_{\delta(\omega)}(\omega)$, there is a random variable
$\bar{\delta}\leq \delta$ almost surely such that for all $x$,
$y\in H$,
$$
x-y\in\mathcal{C}_{\delta(\omega)}(\omega),
$$
implies
$$
\varphi(t,\omega)x-\varphi(t,\omega)y\in
\mathcal{C}_{\bar{\delta}(\theta_t\omega)}(\theta_t\omega), \;\;
{\rm for \;\; almost \;\; all\;\;}\omega\in\Omega.
$$
Then the random dynamical system $\varphi(t,\omega)$ is called have
the cone invariance property for the cone
$\mathcal{C}_{\delta(\omega)}(\omega)$.
\end{defn}

\begin{remark}
Both asymptotic completeness and cone invariance are   important
tools to study the inertial manifold of deterministic infinite
dimensional systems  \cite{Rob93, Rob96, KokSie}. Here we modified
both concepts for   random systems.
\end{remark}

\begin{remark}
Almost sure asymptotic completeness   describes the attracting
property of $\mathcal{M}(\omega)$ for RDS $\varphi(t,\omega)$.
When this property holds, the infinite dimensional system
$\varphi(t,\omega)$ can be reduced to a finite dimensional system
on $\mathcal{M}(\omega)$, and the asymptotic behavior of
$\varphi(t,\omega)$ can be determined by that of the
reduced system on $\mathcal{M}(\omega)$.\\
\end{remark}

 For the random dynamical system $\varphi(t,\omega)$ generated by
the random PDE (\ref{rpde1}), we have the following result.

\begin{lemma} \label{approximation}
For small Lipschitz constant $L_F$ random dynamical system
$\varphi(t,\omega)$ possesses the cone invariance property for a
cone with a deterministic positive constant $\delta$. Moreover if
there exists $t_0>0$ such that for $x$, $y\in H$ and
$$ \varphi(t_0,\omega)x-\varphi(t_0,\omega)y \;\notin
\; \mathcal{C}_\delta(\theta_{t_0}\omega),
$$
then
$$
|\varphi(t,\omega)x- \varphi(t,\omega)y|\leq
D(\omega)e^{-kt}|x-y|,\;\; 0\leq t\leq t_0,
$$
where $D(\omega)$ is a positive tempered random variable and $k=
\beta-L_F-\delta^{-1}L_F>0$.
\end{lemma}
Note that the smallness condition on the Lipschitz constant $L_F$
is specifically defined in (\ref{c}) below.

\begin{proof}
Let $v$, $\bar{v}$ be two solutions of (\ref{rpde1}) and
$p=v_c-\bar{v}_c$, $q=v_s-\bar{v}_s$, then
\begin{equation}\label{p}
\dot{p}= A_cp+z(\theta_t\omega)p+g_c(\theta_t\omega,
v_c+v_s)-g_c(\theta_t\omega, \bar{v}_c+\bar{v}_s),
\end{equation}
\begin{equation}\label{q}
\dot{q}=  A_sq+z(\theta_t\omega)q+g_s(\theta_t\omega,
v_c+v_s)-g_s(\theta_t\omega, \bar{v}_c+\bar{v}_s).
\end{equation}
From (\ref{p}), (\ref{q}) and by the property of $A$ and $F$ we have
\begin{equation}\label{p2}
\frac{1}{2}\frac{d}{dt}|p|^2\geq - \alpha |p|^2
+z(\theta_t\omega)|p|^2-L_F|p|^2-L_F|p|\cdot|q|
\end{equation}
and
\begin{equation}\label{q2}
\frac{1}{2}\frac{d}{dt}|q|^2\leq - \beta |q|^2
+z(\theta_t\omega)|q|^2+L_F|q|^2+L_F|p|\cdot|q|.
\end{equation}
Then (\ref{q2})$-\delta^2\times$ (\ref{p2}), we have
\begin{eqnarray*}
\frac{1}{2}\frac{d}{dt}(|q|^2-\delta^2|p|^2)&\leq &
- \beta|q|^2+z(\theta_t\omega)|q|^2+L_F|q|^2+L_F|p|\cdot|q|+\\
&& \alpha\delta^2|p|^2-z(\theta_t\omega)\delta^2|p|^2+\delta^2
L_F|p|^2+\delta^2 L_F|p|\cdot|q|.
\end{eqnarray*}
Note that if $(p, q)\in \partial \mathcal{C}_\delta(\omega)$ (the
boundary of the cone $\mathcal{C}_\delta(\omega)$), then
$|q|=\delta|p|$ and
\begin{equation*}
\frac{1}{2}\frac{d}{dt}(|q|^2-\delta^2|p|^2)\leq ( \alpha-
\beta+2L_F+\delta L_F+\delta^{-1}L_F)|q|^2.
\end{equation*}
If $L_F$ is small enough such that
\begin{equation}\label{c}
 \alpha- \beta+2L_F+\delta
L_F+\delta^{-1}L_F<0,
\end{equation}
then $|q|^2-\delta^2|p|^2$ is decreasing on $\partial
\mathcal{C}_\delta(\omega)$. Thus it is obvious that  whenever
$x-y \in\mathcal{C}_{\delta}(\omega)$,
$\varphi(t,\omega)x-\varphi(t,\omega)y$ can not leave
$\mathcal{C}_\delta(\theta_t\omega)$.

We now prove the second claim. If there is $t_0>0$ such that
$\varphi(t_0,\omega)x-\varphi(t_0,\omega)y\;\notin\;\mathcal{C}_\delta(\theta_{t_0}\omega)$,
the cone invariance yields
$$
\varphi(t,\omega)x-\varphi(t,\omega)y\; \notin\;
\mathcal{C}_\delta(\theta_t\omega),\;\; 0\leq t\leq t_0,
$$
that is
$$
|q(t)|>\delta|p(t)|,\;\; 0\leq t\leq t_0.
$$
Then by (\ref{q2}) we have
$$
\frac{1}{2}\frac{d}{dt}|q|^2\leq -(
\beta-L_F-\delta^{-1}L_F-z(\theta_t\omega))|q|^2, \;\;0\leq t\leq
t_0.
$$
Hence
$$
|p(t)|^2<\frac{1}{\delta^2}|q(t)|^2\leq
\frac{1}{\delta^2}e^{-2kt+\int_0^tz(\theta_s\omega)ds},\;\;0\leq
t\leq t_0.
$$
Then by the property of $z(\theta_t\omega)$ there is a tempered
random variable $D(\omega)$ such that
$$
|\varphi(t,\omega)x- \varphi(t,\omega)y|\leq
D(\omega)e^{-kt}|x-y|,\;\; 0\leq t\leq t_0.
$$

This completes the proof of the lemma.
\end{proof}

\vskip 0.3cm

Before we prove   Theorem \ref{theorem}, we need the following
lemma, which implies the backward solvability of the system
(\ref{rpde1}) restricted on the invariant manifold
$\widetilde{\mathcal{M}}(\omega)$.

For any given final time $T_f>0$, consider the following system
for $t\in [0, T_f]$
\begin{eqnarray}
\dot{v}_c&=& A_c v_c+z(\theta_t\omega)v_c+ g_c(\theta_t\omega,
v_c+v_s),
\;\; v_c(T_f)=\xi\in H_c \\
\dot{v}_s&=& A_s v_s+z(\theta_t\omega)v_s+g_s(\theta_t\omega,
v_c+v_s),\;\; v_s(0)=h^s(v_c(0))
\end{eqnarray}
where $h^s$ is defined as (\ref{h}). Rewrite the above problem in
the following equivalent integral form
\begin{eqnarray}\label{b1}
v_c(t)=e^{A_c(t-T_f)+\int_{T_f}^tz(\theta_\tau\omega)d\tau}\xi+
\int_{T_f}^te^{A_c(t-\tau)+\int_\tau^tz(\theta_{\varsigma}\omega)d\varsigma}
g_c(\theta_\tau\omega, v(\tau))d\tau,
\end{eqnarray}
\begin{eqnarray}\label{b2}
 v_s(t)=e^{A_st+\int_0^tz(\theta_\tau\omega)d\tau}h^s(v_c(0))
 +\int_0^te^{A_s(t-\tau)+\int_\tau^tz(\theta_\varsigma\omega)d\varsigma}
 g_s(\theta_\tau\omega, v(\tau))d\tau
\end{eqnarray}
$t\in [0,T_f]$.

\begin{lemma}\label{back}
Let (\ref{c1}) hold. Then for any $T_f>0$, (\ref{b1}), (\ref{b2})
has a unique solution $(v_c(\cdot),v_s(\cdot))\in C(0, T_f;
H_c\times H_s)$. Moreover for any $t\geq 0$,
$(v_c(t,\theta_{-t}\omega), v_s(t,   \theta_{-t}\omega))\in
\widetilde{\mathcal{M}}(\omega)$ for almost all $\omega\in\Omega$.
\end{lemma}
\begin{proof}
The existence and uniqueness on small time interval can be
obtained by a contraction  argument, as in Lemma 3.3 of
\cite{DLSch03}.   Then the solution can be extended to any time
interval; see Theorem 3.8 of \cite{DLSch03}.
\end{proof}
\bigskip

Now we complete the proof of the main result of this paper.

\bigskip

 \textbf{Proof of the main result: Theorem \ref{theorem} in \S \ref{s2} }

\bigskip
 It remains only to prove the almost sure  asymptotic completeness of
$\widetilde{\mathcal{M}}(\omega)$.  We fix a $\omega\in\Omega$.
Consider a solution
$$
v(t,\theta_{-t}\omega)=(v_c(t,\theta_{-t}\omega),
v_s(t,\theta_{-t}\omega))
$$
of (\ref{rpde1}). For any $\tau>0$  by Lemma \ref{back} we can find
a solution of (\ref{rpde1}) $\bar{v}(t,\theta_{-t}\omega)$, lying on
$\widetilde{\mathcal{M}}(\omega)$ such that
$$
\bar{v}_c(\tau, \theta_{-\tau}\omega)=v_c(\tau,
\theta_{-\tau}\omega).
$$
Then $\bar{v}(t, \theta_{-t}\omega)$ depends on $\tau>0$. Write
\begin{equation*}
\bar{v}_c(0; \tau,\omega):=\bar{v}_c(0,\omega)
\end{equation*}
and
\begin{equation*}
\bar{v}_s(0;\tau,\omega):=\bar{v}_s(0,\omega).
\end{equation*}
By the construction of $\widetilde{\mathcal{M}}(\omega)$
\begin{eqnarray*}
|\bar{v}_s(0;\tau,\omega)|&\leq &\int_0^\infty e^{-\beta
r-\int_0^{-r}z(\theta_\varsigma\omega)d\varsigma}|g_s(\theta_r\omega,
v^*(-r))|dr\\
 &\leq&L_F\int_0^\infty e^{-(\beta+\eta) r}e^{\eta r}
 e^{-\int_0^{-r}z(\theta_\varsigma\omega)d\varsigma}|v^*(-r)|dr\\
 &:=& N_{L_F}(\omega)\\
 &\leq& L_F|v^*|_{C^-_\eta}\int_0^\infty e^{-(\beta+\eta) r}dr. \;\;\;\;\;(\beta+\eta>0 {\rm \;by\; the\; choice\; of\;} \eta)
\end{eqnarray*}
 It is easy to see that $N_{L_F}(\omega)$ is a finite tempered random variable
and  $N_{L_F}(\omega)\sim  O(L_F)$
 almost surely. And since
$\bar{v}_c(\tau,\theta_{-\tau}\omega)=v_c(\tau,\theta_{-\tau}\omega)$,
by the cone invariance
$$
v(t,\theta_{-t}\omega)-\bar{v}(t,\theta_{-t}\omega)\notin
\mathcal{C}_{\delta}(\omega),\;\; 0\leq t\leq \tau.
$$

Let $S(\omega)=\{\bar{v}_c(0; \tau, \omega):\tau >0 \}$. Notice that
\begin{eqnarray*}
|\bar{v}_c(0; \tau,\omega)-v_c(0,
\omega)|&<&\frac{1}{\delta}|\bar{v}_s(0; \tau, \omega)-v_s(0,\omega)|\\
&\leq& \frac{1}{\delta}\big(N_{
 L_F}(\omega)+|v_s(0,\omega)|\big).
\end{eqnarray*}
Then $S$ is a random bounded set in finite dimensional space, that
is for almost all $\omega\in\Omega$, $ S(\omega)$ is a bounded set
in $\R^n$ and the bound may not be uniform to $\omega\in\Omega$. But
for almost all $\omega\in\Omega$ we can pick out a sequence
$\tau_m\rightarrow \infty $ such that
$$
\lim_{m\rightarrow \infty}\bar{v}_c(0; \tau_m,\omega)=V_c(\omega).
$$
Moreover $V(\omega)$ is measurable with respect to $\omega$. Define
 $V(t,\theta_{-t}\omega)=(V_c(t,\theta_{-t}\omega),$\\$
V_s(t,\theta_{-t}\omega))$ be a solution of (\ref{rpde1}) with
$V(0,\omega)=(V_c(\omega), h^s(V_c(\omega),\omega))$. Then $V(t,
\theta_{-t}\omega)\in\widetilde{\mathcal{M}}(\omega)$ and it is
easy to check by a contradiction argument that
$$
v(t,\theta_{-t}\omega)-V(t,\theta_{-t}\omega)
\notin\mathcal{C}_\delta(\omega), \;\;0\leq t <\infty
$$
which means the almost sure asymptotic completeness of
$\widetilde{\mathcal{M}}(\omega)$.

 This finishes the proof of Theorem \ref{theorem}.

\medskip

\begin{remark}
This theorem implies that the random system (\ref{rpde1}) is
asymptotically reduced to the following random ordinary differential
equation on $H_c$
\begin{equation}\label{detequ}
\dot{v}_c=  A_c v_c+z(\theta_t\omega)v_c+g_c(\theta_t\omega,
v_c+h^s(\omega, v_c))
\end{equation}
where $h^s$ is given by the equation (\ref{h}). Thus the stochastic
PDE     (\ref{SEE}) is asymptotically reduced to the following
finite dimensional non-autonomous stochastic differential system on
the invariant manifold $H_c$:
\begin{equation}\label{detesde}
du_c=  [A_c u_c+F_c(u_c+\bar{h}^s(u_c))]dt+u_c\circ dW(t)
\end{equation}
where $\bar{h}^s(u_c)=e^{z(\theta_t\omega)}h^s(\theta_t\omega, v_c
)$.
\end{remark}

\begin{remark} \label{true}
From the proof of the Theorem  \ref{theorem}, we see that the
dynamical approximation estimate holds for the random partial
differential equation (\ref{rpde1}). In other words, we may apply
Theorem \ref{theorem} to random partial differential equations
like (\ref{rpde1}), as long as appropriate conditions are
satisfied.

\end{remark}


\section{Applications: Detecting  stationary states} \label{s5}

In this final section, as an application of Theorem \ref{theorem}
and its consequences, we consider stationary solutions for a
stochastic hyperbolic equation.

We intend to detecting the stationary states of the following
hyperbolic equation driven by multiplicative white noise,
\begin{equation}\label{HypM}
 u_{tt}+a u_t=\nu\Delta u+bu+f(u)+  u\circ\dot{W},\;
\; {\rm in\;\; I}
\end{equation}
with
$$
u(0)=u_0,\;\;u_t(0)=u_1,\;\; u(0)=u(2\pi)=0
$$
where ${\rm I}$ is taken as the interval $(0, 2\pi)$ for
simplicity,    both $\nu$ and $ a>0$ are positive constants, and
$W(t)$ is  a scalar Wiener process. Moreover, $f\in
C^{1,1}(\R,\R)$ is bounded with global Lipschitz constant $L_f$.
For example $f(x)=\sin x$, which yields the Sine-Gordon equation.

 We study the existence of the
stationary solutions of (\ref{HypM}) by reducing  the system to a
finite dimensional system on an invariant manifold, which is
almost surely asymptotically complete. In fact  by Theorem
\ref{theorem} in \S \ref{s2}, the dynamical behavior of
(\ref{HypM}) restricted on the invariant manifold, is determined
by a stochastic ordinary differential system.   The existence and
stability of the stationary solutions of the stochastic hyperbolic
equation can be obtained from that of the stochastic ordinary
differential system. A similar result for parabolic system with
large diffusivity holds   \cite{Hale86}.

Let $\mathcal{H}=H_0^1({\rm I})\times L^2({\rm I})$. Rewrite the
system (\ref{HypM}) as the following one order stochastic
evolutionary equation in $\mathcal{H}$
\begin{eqnarray}
du&=&vdt\label{shyp1}\\
dv&=&\big[-\nu A u+bu-a v+f(u)\big]dt+  u\circ dW(t) \label{shyp2}
\end{eqnarray}
where $A=-\Delta$ with Dirichlet boundary condition on ${\rm I}$.
Note that Theorem \ref{theorem} can not be applied directly to the
system (\ref{shyp1})-(\ref{shyp2}). However, by Remark \ref{true},
we
 can still have the same result of the Theorem \ref{theorem} for the
stochastic hyperbolic system (\ref{HypM}). This  will be clear
 after we transform (\ref{shyp1})-(\ref{shyp2}) into the form of
(\ref{rpde1}).

First we prove that the system (\ref{shyp1})-(\ref{shyp2}) generates
a continuous random dynamical system in $\mathcal{H}$. Let
$\phi_1(t)=u(t)$, $\phi_2(t)=u_t(t)-  u(t) z(\theta_t\omega)$ where
$z(\omega)$ is the stationary solution of (\ref{OU}). Then we have
the following random evolutionary equation
\begin{eqnarray}
d\phi_1&=&\big[\phi_2+  \phi_1 z(\theta_t\omega)\big]dt\label{rhpy1}\\
d\phi_2&=&\big[\nu\D\phi_1+b\phi_1-a\phi_2+f(\phi_1)\big]dt+\big[(
z(\theta_t\omega)-  a z(\theta_t\omega)-\nonumber\\
&& z^2(\theta_t\omega))\phi_1-  z(\theta_t\omega)\phi_2
\big]dt\label{rhyp2}.
\end{eqnarray}
By a standard Galerkin approximation procedure as in \cite{Tem97},
the system (\ref{rhpy1})-(\ref{rhyp2}) is wellposed. In fact we
give a prior estimates. Multiplying (\ref{rhpy1}) with $\nu\phi_1$
in $H_0^1({\rm I})$ and (\ref{rhyp2}) with $\phi_2$ in $L^2(\rm
I)$. Since $f$ is bounded by a simple calculation we have
\begin{equation*}
\frac{d}{dt}\big[\nu|\phi_1|^2_{H_0^1({\rm
I})}+|\phi_2|^2_{L^2(\mathcal{D})}\big]\leq
C\big(1+|z(\theta_t\omega)|^2+|z(\theta_t\omega)|\big)\big[\nu|\phi_1|^2_{H_0^1({\rm
I})}+|\phi_2|^2_{L^2({\rm I})}\big]
\end{equation*}
for appropriate constant $C$. Then for any $T>0$, $(\phi_1, \phi_2)$
is bounded in $L^\infty(0, T; \mathcal{H})$ which ensures the
weak-star convergence by the Lipschitz property of $f$, for the
detail see \cite{KSch99}. Let $\Phi(t,\omega)=(\phi_1(t,\omega),
\phi_2(t,\omega))$, then $\varphi(t,\omega, \Phi(0))=\Phi(t,\omega)$
defines a continuous random dynamical system in $\mathcal{H}$.
Notice that the stochastic system (\ref{shyp1})-(\ref{shyp2}) is
conjugated to the random system (\ref{rhpy1})-(\ref{rhyp2}) by the
homeomorphism
$$
T(\omega, (u, v))=(u, v+uz(\omega)),\;\; (u, v)\in\mathcal{H}
$$
with inverse
$$
T^{-1}(\omega, (u, v))=(u, v-uz(\omega)),\;\; (u,
v)\in\mathcal{H}.
$$
Then $\hat{\varphi}(t,\omega, (u_0, v_0))=T(\theta_t\omega,
\varphi(t,\omega, T^{-1}(\omega,(u_0, v_0) )))$ is the random
dynamical system generated by (\ref{shyp1})-(\ref{shyp2}). For
more relation about the two systems see \cite{KSch99}.

Define $$ \mathcal{A}_\nu=\left(
\begin{array}{cc}
    0, &  -id_{L^2({\rm I})}  \\
    \nu A-b, &  a\\
\end{array}
\right),\;\; F(\Phi)= \left(
\begin{array}{cc}
    0   \\
    f(\phi_1)  \\
\end{array}
\right),$$
$$ Z(\theta_t\omega)=\left(%
\begin{array}{cc}
       z(\theta_t\omega), & 0 \\

z(\theta_t\omega)-  a z(\theta_t\omega)-  z^2(\theta_t\omega), & -
z(\theta_t\omega)
\end{array}
\right),
$$
where $-id_{L^2({\rm I})}$ is the identity operator on the Hilbert
space $L^2({\rm I})$.

 Then (\ref{rhpy1})-(\ref{rhyp2}) can be written as
\begin{equation}\label{RHPY}
\frac{d\Phi}{dt}=-\mathcal{A}_\nu\Phi+Z(\theta_t\omega)\Phi+F(\Phi)
\end{equation}
which is in the form of (\ref{rpde1}). Thus by Remark \ref{true},
we can still apply Theorem  \ref{theorem} here.

The eigenvalues of the operator $A$ are $\lambda_k=k^2$
 with corresponding eigenvectors $\sin kx$, $k=1,2, \cdots$. Then
the operator $\mathcal{A}_\nu$ has the eigenvalues
$$
\delta_k^\pm=\frac{a}{2}\pm\sqrt{\frac{a^2}{4}+b-\nu
k^2},\;\;k=1,2,\cdots
$$
and corresponding eigenvectors are $(1, \delta_k^\pm)\sin kx $.
Define subspace of $\mathcal{H}$
$$
\mathcal{H}_{c}= {\rm span} \{(1,0)^T\sin 2x\},
\;\;\mathcal{H}_{u}= {\rm span} \{(1, \delta_1^-)^T\sin x\}
$$
and $\mathcal{H}_{cu}=\mathcal{H}_{c}\oplus \mathcal{H}_{u}$.
Write the projections from $\mathcal{H}$ to $\mathcal{H}_{c}$,
$\mathcal{H}_{u}$ and $\mathcal{H}_{cu}$ as $P_c$, $P_u$ and
$P_{cu}$ respectively. We also use the subspaces
$$
\mathcal{H}^-_{c}= {\rm span} \{(1,\delta_2^\pm)^T\sin 2x\},\;\;
\mathcal{H}^-_{u}= {\rm span} \{(1,\delta_1^\pm)^T\sin 2x\}
$$
with the projections $P^-_c$ and $P_u^-$ from $\mathcal{H}$ to
$\mathcal{H}^-_{c}$ and $\mathcal{H}^-_{u}$ respectively. Let
$\mathcal{H}^-_{cu}=\mathcal{H}^-_{c}\oplus \mathcal{H}^-_{u}$ and
$\mathcal{H}^-_s={\rm span}\{ (1, \delta_k^\pm)^T \sin kx,
k\in\Z^+\setminus \{1,2\}\}$.

 Here we consider a special case that $b=4\nu$ and
$\nu=a^2/4$. Then the operator $\mathcal{A}_\nu$ has one zero
eigenvalue $\delta_2^-=0$, one negative eigenvalue
$\delta_1^-=\frac{a}{2}$ and the others are all complex numbers
with positive real part. Since $\{(1, \delta_k^\pm)\sin kx \}$ are
not orthogonal, we introduce a new inner product (see
\cite{Mora87,Qin01,SMV87}), which defines an equivalent norm on
$\mathcal{H}$. For $(u_1, v_1)$, $(u_2,v_2 )\in
\mathcal{H}^-_{cu}$, define
$$
\langle (u_1, v_1), (u_2,v_2 ) \rangle_{cu}=2\big(\langle-Au_1,
u_2\rangle_0+\frac{a^2}{4}\langle u_1, u_2\rangle_0+\langle
 \frac{a}{2}u_1+v_1, \frac{a}{2}u_2+v_2 \rangle_0
 \big)
$$
and for $(u_1, v_1)$, $(u_2,v_2 )\in \mathcal{H}^-_{s}$ define
$$
\langle (u_1, v_1),(u_2,v_2 ) \rangle_{s}=2\big(\langle Au_1,
u_2\rangle_0-\frac{a^2}{4}\langle u_1, u_2\rangle_0+\langle
 \frac{a}{2}u_1+v_1, \frac{a}{2}u_2+v_2 \rangle_0
 \big)
$$
where $\langle\cdot, \cdot\rangle_0$ is usual inner product in
$L^2({\rm I})$. Then we introduce the following new inner product
in $\mathcal{H}$ defined by
$$
\langle U, \tilde{U}\rangle_{\mathcal{H}}=\langle
U_{cu},\tilde{U}_{cu}
   \rangle_{cu}+\langle U_{s},\tilde{U}_{s}   \rangle_s
$$
for $U=U_{cu}+U_s$, $\tilde{U}=\tilde{U}_{cu}+\tilde{U}_s\in
\mathcal{H}$. And the new norm
$$
\|U\|^2_\mathcal{H}=\langle U, U\rangle_{\mathcal{H}}
$$
is equivalent to the usual norm $\|\cdot\|_{H_0^1\times L^2(D)}$.
Moreover $\mathcal{H}_{cu}$ and
$\mathcal{H}_{s}=\mathcal{H}\ominus\mathcal{H}_{cu}$ is
orthogonal. The new norm has the following properties\cite{Qin01}:

\begin{enumerate}
    \item $\|U\|_\mathcal{H}=\sqrt{2}\|v\|_{L^2(D)}$, for $U=(0,
    v)\in\mathcal{H}$.
    \item $\|U\|_\mathcal{H}\geq \frac{\sqrt{2}a}{2}\|u\|_{L^2(D)}$, for $U=(u,
    v)\in\mathcal{H}$.
    \item In terms of the new norm the Lipschitz constant of $f$ is
    $\frac{2L_f}{a}$.
\end{enumerate}

A   calculation yields that, in terms of the new norm, the
semigroup $\mathcal{S}(t)$ generated by $\mathcal{A}_\nu$
satisfies (\ref{expdic1})-(\ref{expdic2}) with $\alpha=0$,
$\beta=\frac{a}{2}$ and $\dim H_c=\dim \mathcal{H}_{cu}=2$. Taking
$\eta=\frac{-\beta-\alpha}{2}$, condition (\ref{c1}) for
(\ref{HypM}) becomes
\begin{equation}\label{c11}
\frac{8L_f}{a^2}<1
\end{equation}
which holds if $a$ is large enough.

Define the space $C_\eta^-$ with $z(\theta_t\omega)$ replaced by
$Z(\theta_t\omega)$. Then by the random invariant manifold theory in
\S 3 above, we know the system (\ref{RHPY}) has  a two dimensional
random invariant manifold $\widetilde{\mathcal{M}}(\omega)$ in
$C_\eta^-$ provided $a$ is large enough. And
$\widetilde{\mathcal{M}}(\omega)$ can be denoted as the graph of a
random Lipschitz map $h^s:\mathcal{H}_{cu}\rightarrow \mathcal{H}_s$
with Lipschitz constant $CL_f/a$. Then by Theorem \textbf{A}, the
dynamical behavior of (\ref{RHPY}), that is, (\ref{HypM}), is
determined by the following reduced system on $\mathcal{M}(\omega)$
which is a finite dimensional random system
\begin{equation}\label{reduced}
\frac{d\Phi_{cu}}{dt}=-P_{cu}\mathcal{A}_{\nu}\Phi_{cu}+Z(\theta_t\omega)\Phi_{cu}+P_{cu}F(\Phi_{cu}+h^s(\Phi_{cu})).
\end{equation}
provided $a$ is large enough. To see the reduced system more clear,
projecting the system (\ref{HypM}) to span$\{\sin x \}$ and
span$\{\sin 2x\}$ by the projection $\Pi_1$ and $\Pi_2$ respectively
we have the following system
\begin{eqnarray}
\dot{u}_1&=&v_1 \label{u1} \\
\dot{v}_1&=&\frac{3a^2}{4}u_1-a v_1+\Pi_1 f(u)+
u_1\circ \dot{W}(t)  \label{v1}\\
\dot{u}_2&=&v_2\label{u2}\\
\dot{v}_2&=&-a v_2+\Pi_2 f(u)+  u_2\circ \dot{W}(t).\label{v2}
\end{eqnarray}
By the new inner product in $\mathcal{H}$, define the following
new inner product in $\R^2$ as
$$
\big((x, y), (\tilde{x},\tilde{y})\big)=2\big[\frac{a^2}{4} x
\tilde{y}+ (\frac{a}{2}x+y)(\frac{a}{2}\tilde{x}+\tilde{y})\big].
$$
Then projecting the above system (\ref{u1})-(\ref{v2}) to
$\mathcal{H}_{cu}$ by the new inner product in $\R^2$, we have the
reduced system
\begin{eqnarray}
\dot{u}_1&=&\frac{a}{2}u_1+\frac{2\sqrt{10}}{5}\Pi_1
f(u_1+u_2+\bar{h}^s(u_1, u_2))+
u_1\circ \dot{W}(t)\label{u0} \\
 \dot{u}_2&=&\Pi_2 f(u_1+u_2+\bar{h}^s(u_1, u_2))+
u_2\circ \dot{W}(t)\label{v0}
\end{eqnarray}
where $\bar{h}^s=e^{z(\omega)}h^s$.

Then the dynamical behavior of system (\ref{HypM}) is determined by
that of the system (\ref{u0})-(\ref{v0}) if $a$ is large. So if
system (\ref{u0})-(\ref{v0}) has a stationary solution, so does
system (\ref{HypM}). And the stability property of the stationary
solution is also determined by that of (\ref{u0})-(\ref{v0}).







\end{document}